\documentclass[11pt]{article}
\usepackage[all]{xy}
\usepackage{amsfonts,amsmath,oldgerm,amssymb,amscd}
\newcommand{\ra}{\rightarrow}

\newcommand{\by}[1]{\stackrel{#1}{\ra}}

\newcommand{\surj}{\ra\!\!\!\ra}	\newcommand{\inj}{\hookrightarrow}

\newcommand{\ol}{\overline}		\newcommand{\wt}{\widetilde}
\newcommand{\iso}{\by \sim}

\newtheorem{theorem}{Theorem}[section]
\newtheorem{proposition}[theorem]{Proposition}
\newtheorem{lemma}[theorem]{Lemma}

\newtheorem{corollary}[theorem]{Corollary}
\newtheorem{question}[theorem]{Question}

\newcommand{\ga}{\alpha}		\newcommand{\gb}{\beta}
\newcommand{\gc}{\chi}			\newcommand{\gd}{\delta}
\newcommand{\gf}{\varphi}		
\newcommand{\gj}{\square}		\newcommand{\gl}{\lambda}

\newcommand{\gt}{\theta}		
		\newcommand{\gw}{\wedge}

		\newcommand{\gT}{\Theta}

\newcommand{\BQ}{\mbox{$\mathbb Q$}}

	\newcommand{\CJ}{\mbox{$\mathcal J$}}

	\newcommand{\CR}{\mbox{$\mathcal R$}}

\newcommand{\op}{\mbox{$\oplus$}}	
	\newcommand{\Spec}{\mbox{\rm Spec\,}}
\newcommand{\hh}{\mbox{\rm ht\;}}	\newcommand{\Hom}{\mbox{\rm Hom}}
  	
\newcommand{\Um}{\mbox{\rm Um}}		\newcommand{\SL}{\mbox{\rm SL}}
\newcommand{\GL}{\mbox{\rm GL}}		\newcommand{\ot}{\mbox{$\otimes$}}
\newcommand{\Aut}{\mbox{\rm Aut}}

\def\notation{\refstepcounter{theorem}\paragraph{Notation \thetheorem}}

\begin{document}

\begin{center}
{\bf \Large Euler class group of a Laurent polynomial ring : local case}
\vspace{.2in} \\
	{\large Manoj  Kumar Keshari} 
\vspace{.1in}\\
{\small 
Department of Mathematics, IIT Mumbai, Mumbai - 400076, India;\;
	keshari@math.iitb.ac.in} 
\end{center}


\section{Introduction}

Let $A$ be a commutative Noetherian ring of dimension $d$. A classical
result of Serre \cite{S1} asserts that if $P$ is a
projective $A$-module of rank $>d$, then $P$ has a unimodular
element. It is well known that this result is not true in general if
rank $P=d=\dim A$. Therefore, it is interesting to know the obstruction for
projective $A$-modules of rank $=\dim A$ to have a unimodular element. 

Let $A$ be a commutative Noetherian ring of dimension $n$ containing
$\BQ$ and let $P$ be a projective $A$-module of rank $n$.  In
\cite{BR3}, an abelian group $E(A)$, called the {\it Euler class
group} of $A$ is defined and it is shown that $P$ has a unimodular
element if and only if the Euler class of $P$ in $E(A)$ vanishes (see
\cite{BR3} for the definition of Euler class of $P$).

In view of the above result \cite{BR3}, we can ask the following:

\begin{question}
Let $A$ be a commutative Noetherian ring containing $\BQ$.  Let $P$ be
a projective $A$-modules of rank $r < \dim A$ having trivial
determinant. What is the obstruction for $P$ to have a unimodular
element?
\end{question}

Let $R$ be a commutative Noetherian ring of dimension $n$ containing $\BQ$.
In \cite{Das}, an abelian group $E(R[T])$, 
called the Euler class group of $R[T]$ is defined and it is shown that if $P$ 
is a projective $R[T]$-module of rank $n=\dim R[T] -1$ 
with trivial determinant,
then $P$ has a unimodular element
if and only if the Euler class of $P$ in $E(R[T])$ vanishes, thus answering the
above question in the case $r=\dim A-1$ and $A=R[T]$.

In this paper, we prove results similar to \cite{Das} for
the ring $R[T,T^{-1}]$ under the assumption that height
of the Jacobson radical of $R$ is $\geq 2$.
 More precisely, we define the Euler
class group of $R[T,T^{-1}]$ and prove that 
if $\wt P$ is a projective $R[T,T^{-1}]$-module of rank $n=\dim R$ 
with trivial determinant, then $\wt P$ has a unimodular element
if and only if the Euler class of $\wt P$ in $E(R[T,T^{-1}])$ vanishes 
(\ref{cor-1}). 

In the appendix, we prove the following ``Symplectic'' cancellation
theorem (\ref{88}) (it is used in Section $7$)
which is a generalization of (\cite{Bh-2}, Theorem 4.8), where
it is proved in the polynomial ring case. 

\begin{theorem}
Let $B$ be a ring of dimension $d$ and $A=B[Y_1,\ldots,Y_{s},X_1^{\pm
1},\ldots, X_r^{\pm 1}]$. Let $(P,\langle,\rangle)$ be a symplectic
$A$-module of rank $2n > 0.$ If $2n \geq d$, then $ESp(A^2 \perp P,
\langle,\rangle)$ acts transitively on $\Um(A^2 \oplus P).$
\end{theorem}

As an application, we get the following result (\ref{87}), which
gives a partial answer to a question of Weibel (\cite{weibel}, Introduction).

\begin{theorem}
Let $R$ be a ring of dimension $2$
and $A=R[X_1,\ldots,X_r,Y_1^{\pm 1},\ldots,Y_{s}^{\pm 1}]$.
Assume $A^2$ is cancellative. Then every projective 
$A$-modules of rank $2$ with trivial determinant is cancellative.
\end{theorem} 

\section{Preliminaries}

All the rings considered in this paper are assumed to be
commutative Noetherian and  all the modules are finitely generated. We
denote the Jacobson radical of $A$ by $\CJ(A)$.

Let $B$ be a ring and let $P$ be a projective $B$-module. Recall
that $p\in P$ is called a unimodular element if there exists an $\psi \in
P^*=\Hom_B(P,B)$ such that $\psi(p)=1$. We denote by $\Um(P)$, the set of
all unimodular elements of $P$. 
 
Given an element $\gf\in P^\ast$ and an element $p\in P$, we define an
endomorphism $\gf_p$ as the composite $P\by \gf B\by p P$. 
If $\gf(p)=0$, then ${\gf_p}^2=0$ and hence $1+\gf_p$ is a uni-potent 
automorphism of $P$.

By a {\it transvection}, we mean an automorphism of $P$ of the form
$1+\gf_p$, where $\gf(p)=0$ and either $\gf$ is unimodular in $P^\ast$
or $p$ is unimodular in $P$. We denote by $E(P)$, the subgroup of
$\Aut(P)$ generated by all transvections of $P$. Note that $E(P)$ is a
normal subgroup of $\Aut(P)$.

An existence of a transvection of $P$ pre-supposes that $P$ has a
unimodular element. Now, let $P = B\op Q$, $q\in Q, \alpha\in
Q^*$. Then $\Delta_q(b,q')=(b,q'+bq)$ and
$\Gamma_\alpha(b,q')=(b+\alpha(q'),q')$ are transvections of
$P$. Conversely, any transvection $\Theta$ of $P$ gives rise to a
decomposition  $P=B\op Q$ in such a way that $\Theta = \Delta_q$ or
$\Theta = \Gamma_\alpha$.

We begin by stating two classical results of Serre \cite{S1} and Bass
\cite{Bass} respectively.

\begin{theorem}\label{serre}
Let $A$ be a ring of dimension $d$. Then any projective
$A$-module $P$ of rank $> d$ has a unimodular element. In particular,
if $\dim A=1$, then any projective $A$-module of trivial determinant is
free.
\end{theorem}

\begin{theorem}\label{B}
Let $A$ be a ring of dimension $d$ and let $P$ be a projective
$A$-module of rank $> d$. Then $E(A\op P)$ acts transitively on
$\Um(A\op P)$. In particular, $P$ is cancellative.
\end{theorem}

The following result is due to Lindel (\cite{L}, Theorem 2.6).

\begin{theorem}\label{lindel}
Let $A$ be a ring of dimension $d$ and 
$R = A[T_1,\ldots,T_n,Y_1^{\pm 1},\ldots,Y_r^{\pm 1}]$.
Let $P$ be a projective $R$-module of
rank $\geq$ max $(2,d+1)$. Then $E(P\oplus R)$ acts transitively on
$\Um(P\oplus R)$. In particular, projective
$R$-modules of rank $> d$ are cancellative.
\end{theorem}

The following result is due to Bhatwadekar and Roy (\cite{Bh-Roy},
Proposition 4.1) and is about lifting an automorphism of a projective
module. 

\begin{proposition}\label{trans}
Let $A$ be a ring and $J\subset A$  an ideal. Let $P$ be a
projective $A$-module of rank $n$. 
Then any transvection $\wt \Theta$ of $P/JP$, i.e. $\wt \Theta \in E(P/JP)$,
can be lifted to a (uni-potent) automorphism $\Theta$ of $P$. In
particular, if $P/JP$ is free of rank $n$, then any element $\ol \Psi$ of
$E((A/J)^n)$ can be lifted to $\Psi \in \Aut(P)$. If, in
addition, the natural map $\Um(P) \ra \Um(P/JP)$ is surjective, then
the natural map $E(P)\ra E(P/JP)$ is surjective.
\end{proposition}

The following result is a consequence of a theorem of Eisenbud-Evans
as stated in (\cite{P}, p. 1420).

\begin{lemma}\label{EE}
Let $R$ be a ring and let $P$ be a projective $R$-module of rank
$r$. Let $(\alpha,a) \in (P^\ast \oplus R)$. Then there exists an
element $\beta \in P^\ast$ such that $\hh I_a \geq r$, where
$I=(\alpha+a \gb)(P)$. In particular, if the ideal $(\ga(P),a)$ has
height $\geq r$, then $\hh I \geq r$. Further, if $(\ga(P),a)$ is an
ideal of height $\geq r$ and $I$ is a proper ideal of $R$, then $\hh I
= r$.
\end{lemma}

The following result is  due to Bhatwadekar and Keshari
(\cite{Bh-Manoj}, Lemma 4.4).

\begin{lemma}\label{bm1}
Let $C$ be a ring with $\dim C/\CJ(C) = r$ and
let $P$ be a projective $C$-module of rank $m\geq r+1$. 
Let $I$ and $L$ be ideals of $C$ such that $L\subset I^2$.
Let $ \phi : P \surj I/L$ be a surjection. Then $ \phi$ can be
lifted to a surjection $\Psi : P \surj I$. 
\end{lemma}

The following result is due to Mandal and Raja
Sridharan (\cite{Mandal-Raja}, Theorem 2.3). 

\begin{theorem}\label{Mandal-2}
Let $A$ be a ring and let $I_1,I_2$ be two comaximal ideals
of $A[T]$ such that $I_1$ contains a monic polynomial and $I_2=I_2(0)A[T]$ is
an extended ideal.  Let $I=I_1\cap I_2$. 
Suppose $P$ is a projective $A$-module of rank
$n\geq \dim A[T]/I_1 +2$.  Let $\ga
: P \surj I(0)$ and $\phi : P[T]/I_1P[T] \surj I_1/{I_1}^2$ be two
surjections such that $\phi(0) =\ga \ot A/I_1(0)$. Then there exists a
surjective map $\Psi : P[T] \surj I$ such that $\Psi(0)=\ga$.
\end{theorem}

Now, we state the Addition and Subtraction principles respectively for
arbitrary ring $B$ (\cite{Bh-Manoj}, Theorem 5.6 and Theorem 3.7
respectively). Note that, the following results are valid in the case
$d=n=2$ also (\cite{BR3}, Theorem 3.2 and Theorem 3.3
respectively). 

\begin{proposition}\label{addition}
Let $B$ be a  ring of dimension $d$ and let $I_1,I_2\subset
B$ be two comaximal ideals of height $n$, where $2n \geq d+3$. Let $P=
P_1\op B$ be a projective $B$-module of rank $n$. Let $\Phi : P \surj
I_1$ and $\Psi : P \surj I_2$ be two surjections. Then there exists a
surjection $\Delta : P \surj I_1\cap I_2$ with $\Delta \ot B/I_1 =
\Phi \ot B/I_1$ and $\Delta \ot B/I_2 = \Psi \ot B/I_2$.
\end{proposition}

\begin{proposition}\label{subtract}
Let $B$ be a ring of dimension $d$ and let $I_1,I_2\subset
B$ be two comaximal ideals of height $n$, where $2n \geq d+3$. Let $P=
P_1\op B$ be a projective $B$-module of rank $n$. Let $\Phi : P \surj
I_1$ and $\Psi : P \surj I_1 \cap I_2$ be two surjections such that
$\Phi \ot B/I_1 = \Psi \ot B/I_1$. Then there exists a
surjection $\Delta : P \surj I_2$ such that  $\Delta \ot B/I_2 =
\Psi \ot B/I_2$.
\end{proposition}

We end this section by recalling some results from (\cite{BR3} 4.2,
4.3, 4.4) for later use.   

\begin{theorem}\label{bha}
Let $B$ be a ring of dimension $n\geq 2$ containing
$\BQ$. Let $J$ be an ideal of $B$ of height $n$ such that $J/J^2$ is
generated by $n$ elements. Let $w_J : (B/J)^n \surj J/J^2$.  Let $P$
be a projective $B$-module of rank $n$ with trivial determinant and
$\chi:B\iso \gw^n P$. Then the following holds:

$(1)$ If $(J,w_J)=0$ in $E(B)$, then $w_J$ can be lifted to a
surjection from $B^n$ to $J$.

$(2)$ Suppose $e(P,\chi)=(J,w_J)$ in $E(B)$. Then there exists a
surjection $\ga : P\surj J$ such that $(J,w_J)$ is obtained from
$(\ga,\gc)$.

$(3)$ $e(P,\gc)=0$ in $E(B)$ if and only if $P$ has a unimodular element.
\end{theorem}

\section{Some addition and subtraction principle}

We begin with the following result which is proved in 
(\cite{Bh-Raja-1}, Lemma 3.6)
in the case $A$ is an affine algebra over a field, $f=T$ and
$R=A[T]$. Since the same proof works in our case
also, we omit the proof.

\begin{lemma}\label{moving}
Let $A$ be a ring of dimension $d$ and $R=A[T,T^{-1}]$. Let $\wt P$ be
a projective $R$-module of rank $n$, where $2n \geq d+3$. Let $I
\subset R$ be an ideal of height $n$. Let $J\subset I\cap A$ be any
ideal of height $\geq d-n+2$ and let $f\in R$ be any element. Assume
that we are given a surjection $\phi : \wt P \surj I/(I^2f)$. Then $\phi$
has a lift $\wt \phi : \wt P \ra I$ such that $\wt \phi(\wt P)=I''$
satisfies the following properties :

$(1)~ I'' +(J^2f)=I$,

$(2)~I'' =I\cap I'$, where $\hh I' \geq n$ and

$(3)~ I'+(J^2f)=R$.
\end{lemma}
 
\notation Let $A$ be a ring and $R=A[T,T^{-1}]$. We say $f(T)\in A[T]$
is a {\it special monic polynomial} if $f(T)$ is a monic polynomial
with $f(0)=1$.  By $\CR$, we denote the ring obtained from $R$ by
inverting all the special monic polynomials of $A[T]$. It is easy to
see that $\dim \CR=\dim A$.

The following result is an analogue of (\cite{Bh-Manoj}, Lemma 4.5)
for $A[T,T^{-1}]$.

\begin{lemma}\label{bm3}
Let $A$ be a ring with $\dim A/\CJ(A)=r$ and $R=A[T,T^{-1}]$. Let $I$
and $L$ be ideals of $R$ such that $L\subset I^2$ and $L$ contains a
special monic polynomial. Let $Q$ be a projective $R$-module of rank
$m \geq r+1$. Let $\phi: Q\op R \surj I/L$ be a surjection. Then we
can lift $\phi$ to a surjection $\Phi : Q\op R \surj I$ with
$\Phi(0,1)$ a special monic polynomial.
\end{lemma}

\begin{proof}
Let $\Phi'=(\Theta,g)$ be a lift of $\phi$. Let $f\in L$ be a
special monic polynomial. By adding some multiple of $f$
to $g$, we can assume that the lift $\Phi'=(\Theta,g)$
of $\phi$ is such that $g$ is a special monic polynomial.
 Let $C=R/(g)$. Since $A \inj C$ is an integral extension, we
have $\CJ(A) = \CJ(C) \cap A$ and, hence, $A/\CJ(A) \inj C/\CJ(C)
$ is also an integral extension. Therefore, $\dim C/\CJ(C) =r$.
 
Let ``bar'' denote reduction modulo $(g)$. Then $\Theta$ induces a
surjection $\alpha : \ol {Q}  \surj \ol I/\ol L$, which by 
(\ref{bm1}), 
can be lifted to a surjection from $\ol {Q}$ to $\ol I$. 
Therefore, there exists a map $\Gamma : Q \ra I$ such that
$\Gamma(Q)+ (g)=I$ and $(\Theta-\Gamma)(Q)=K \subset L+(g)$.
Hence $\Theta - \Gamma \in K{Q}^*$. This shows that $\Theta -
\Gamma = \Theta_1 +g \Gamma_1$, where $\Theta_1 \in L{Q}^*$ and
$\Gamma_1 \in {Q}^*$. 

Let $\Phi_1 = \Gamma + g \Gamma_1$ and let $\Phi = (\Phi_1,g)$.
Then $\Phi(Q\op R) = \Phi_1(Q)+(g) = \Gamma(Q)+(g) = I$.
Thus, $\Phi : Q \op R \surj I$ is a surjection. Moreover, 
 $\Phi(0,1) = g$ is a special monic polynomial. 
Since $\Phi - \Phi' =
(\Phi_1 - \Theta,0)$, $\Phi_1 -\Theta \in L{Q}^*$ and $\Phi'$ is a
lift of $\phi$, we see that $\Phi$ is a (surjective) lift of $\phi$.
This proves the result. 
$\hfill \gj$
\end{proof}

The proof of the following result is same as of (\cite{Bh-Manoj},
Lemma 4.6) using (\ref{lindel}, \ref{bm3}).
Hence, we omit the proof. 

\begin{lemma}\label{M2}
Let $A$ be a ring of dimension $d$ and $R=A[T,T^{-1}]$. Let $n$ be an
integer such that $2n \geq d+3$. Let $I$ be an ideal of
$R$ of height $n$ such that $I+\CJ(A) R= R$. Assume that
$\hh \CJ(A)\geq d-n+2$. Let $P=Q\op R^2$ be a projective
$R$-module of rank $n$ and let $\phi : P \surj I/I^2$ be a
surjection. If the surjection $\phi\ot \CR : P\ot \CR \surj
I\CR/I^2 \CR$ can be lifted to a surjection from $P\ot \CR$ to
$I \CR$, then $\phi$ can be lifted to a surjection $\Phi : P
\surj I$.
\end{lemma}

\begin{proposition}\label{addition-1}
{\rm (Addition Principle)} Let $A$ be a ring of dimension $d$ and
$R=A[T,T^{-1}]$. Let $I_1,I_2\subset R$ be two comaximal ideals of
height $n$, where $2n \geq d+3$. Let $P=P'\op R^2$ be a projective
$R$-module of rank $n$. Assume that $\hh \CJ(A) \geq d-n+2$. Let $\Phi
: P \surj I_1$ and $\Psi : P \surj I_2$ be two surjections. Then
there exists a surjection $\Delta : P \surj I_1\cap I_2$ with $\Delta
\ot R/I_1 = \Phi \ot R/I_1$ and $\Delta \ot R/I_2 = \Psi \ot R/I_2$.
\end{proposition}

\begin{remark}
Since $\dim R=d+1$, if $2n \geq d+4$, then we can appeal to
(\ref{addition}) for the proof (without the assumption $\hh
\CJ(A) \geq d-n+2$). So, we need to prove the
result only in the case $2n=d+3$. However, the proof given below works
equally well for $2n> d+3$ and hence, allows us to give a unified
treatment. The same remark is also applicable to 
(\ref{subtract-1}). 
\end{remark}

\begin{proof}
{\bf Step 1 :}
Write $I = I_1\cap I_2$. Let $J=(I\cap A) \cap \CJ(A)$. Since $\hh
(I\cap A) \geq n-1 \geq (d-n+2)$, we have  $\hh J \geq d-n+2$. 
The surjections $\Phi$ and $\Psi$ induces a surjection 
$\Gamma : P \surj I/I^2$ with $\Gamma \ot R/I_1 = \Phi\ot R/I_1$ and
$\Gamma \ot R/I_2 = \Psi \ot R/I_2$. 
It is enough to show that $\Gamma$ has a surjective lift from $P$ to $I$.

Applying (\ref{moving}) with $f=1$, we get a lift $\Gamma_1 \in
\Hom_R(P,I)$ of $\Gamma$ such that the ideal $\Gamma_1(P)=I''$
satisfies the following properties: $(1)~I= I'' +J^2$, $(2)~ I'' =I
\cap K$, where $\hh K \geq n$ and $(3)~ K+J=R$.

Since $\dim \CR=d$, applying 
(\ref{addition}) in the ring $\CR$ for the surjections $\Phi \ot
\CR : P\ot \CR \surj I_1 \CR$ and $\Psi \ot \CR : P\ot \CR \surj
I_2 \CR$, we get a surjective map
$\Delta : P\ot \CR \surj I \CR$ such that 
$\Delta \ot
\CR/I_1 \CR= \Phi \ot \CR/I_1 \CR$ and $\Delta \ot
\CR/I_2 \CR = \Psi \ot \CR/I_2 \CR$. 
It is easy to see, from the very construction of $\Gamma$, that
 $\Delta$ is a lift of $\Gamma \ot \CR$.

We have two surjections 
$\Gamma_1 : P \surj I \cap K$ and $\Delta : P\ot \CR 
\surj I \CR$. Since $\Gamma_1$ is a lift of $\Gamma$, we have
 $\Gamma_1 \ot \CR/I \CR = \Delta \ot \CR/I \CR$. 
Applying (\ref{subtract}) in the ring
$\CR$ for the surjections $\Gamma_1 \ot \CR$ and $\Delta$,
we get a surjection $\Delta_1 : P\ot \CR \surj K 
\CR$ with $\Delta_1 \ot \CR/K \CR = \Gamma_1 \ot \CR/K
\CR$. Since $K$ is comaximal with $J$ and hence with $\CJ(A)$,
applying (\ref{M2}), we get a surjection
$\Delta_2 : P \surj K$ which is a lift of $\Gamma_1 \ot R/K : P \surj
K/K^2$.   \\

{\bf Step 2 :}
We have two surjections $\Gamma_1 : P \surj I \cap K$ and $\Delta_2 : P
\surj K$ with $\Gamma_1 \ot R/K = \Delta_2 \ot R/K$.
Recall that $P=P'\op R^2$, $J= (I\cap A)\cap \CJ(A)$, $K$ is comaximal
with $J$ and $\hh J \geq d-n+2$. Write $P_1 = P'\op R$ and $P=P_1\op R$.

Let ``bar'' denote reduction modulo
$J^2$. Then $\ol R = A/J^2[T,T^{-1}]$ and $\dim A/J \leq
d-(d-n+2)=n-2$. Hence applying 
(\ref{lindel}, \ref{trans}),
we can assume that;
after performing some automorphism of $P_1 \op R$,
$\Delta_2(P_1) = R$ modulo  $J^2$ and $\Delta_2((0,1)) \in J^2$. Assume
that $\Delta_2((0,1))= \gl \in J^2$.
Replacing $\Delta_2$ by $\Delta_2 + \gl \Delta_3$ for some $\Delta_3
\in {P_1}^*$, we can assume, by (\ref{EE}), that $\hh
\Delta_2(P_1)=n-1$. 
Let $\Delta_2(p_1)=1$ modulo $J^2$ for some $p_1\in P_1$.
 Further, replacing $\gl$ by $\gl+\Delta_2(p_1)$,  we
can assume that $\gl=1$ modulo $J^2$.     

Let $K_1$ and $K_2$ be two ideals of $R[Y]$ defined by
$K_1=(\Delta_2(P_1),Y+\gl)$ and $K_2=IR[Y]$.
Then $K_1+K_2 = R[Y]$, since $\Delta_2(P_1)+J=R$ and
$J\subset I$. 
Let $K_3=K_1\cap K_2$. Then we have two surjections 
$\Gamma_1 : P \surj K_3(0)=I \cap K$ and
$\Lambda_1 : P[Y] \surj K_1$ defined by $\Lambda_1 = \Delta_2$  on
$P_1$ and $\Lambda_1((0,1))=Y+ \gl$. Then $\Lambda_1(0) = \Gamma_1$
mod $K_1(0)^2$, as $\Delta_2\ot R/K=\Gamma_1 \ot R/K$.
 Also, note that, since $\hh \Delta_2(P_1)=n-1$ and
$\Delta_2(P_1) + \CJ(A)=R$, $\dim R[Y]/K_1 =\dim R/\Delta_2(P_1) \leq
d-n+1 \leq n-2$. Hence applying (\ref{Mandal-2}), we get a
surjection $\Lambda_2 : P[Y] \surj K_3$ with $\Lambda_2(0)= \Gamma_1$.
Putting $Y=1-\gl $, we get a surjection $\wt \Delta = \Lambda_2(1-\gl)
: P \surj I$ with $\wt \Delta \ot R/I = \Gamma_1 \ot R/I$.

Since $\Gamma_1$ is a lift of $\Gamma : P \surj I/I^2$, we have 
$\wt \Delta \ot R/I = \Gamma \ot R/I$.
This proves the result.
$\hfill \gj$
\end{proof}

\begin{proposition}\label{subtract-1}
{\rm (Subtraction Principle)} Let $A$ be a ring of dimension $d$ and
$R=A[T,T^{-1}]$. Let $I_1,I_2\subset R$ be two comaximal ideals of
height $n$, where $2n \geq d+3$. Let $P=P'\op R^2$ be a projective
$R$-module of rank $n$.  Assume that $\hh \CJ(A) \geq d-n+2$. Let
$\Phi : P \surj I_1\cap I_2$ and $\Psi : P \surj I_1$ be two
surjections with $\Phi \ot R/I_1 = \Psi \ot R/I_1$. Then there exists
a surjection $\Delta : P \surj I_2$ with $\Phi \ot R/I_2 = \Delta \ot
R/I_2$.
\end{proposition}

\begin{proof}
Let $J=(I_2\cap A) \cap
\CJ(A)$. Since $\hh (I_2\cap A) \geq n-1$ and $n-1 \geq d-n+2$, we
have $\hh J \geq d-n+2$. We have a surjection $\phi : P \surj
I_2/{I_2}^2$ induced by $\Phi$. Applying 
(\ref{moving}) with $f=1$, we get a lift $\wt \phi \in \Hom(P,I_2)$ of
$\phi$ such that $\wt \phi(P) = I''$ satisfies the following
properties:
$(1)~ I_2=I''+J^2$,
$(2)~ I'' = I_2\cap K$, where $\hh K \geq n$ and 
$(3)~ K+J^2=R$.

We have two surjections $\Phi : P \surj I_1  \cap I_2$ and $\Psi : P
\surj I_1$ with $\Phi \ot R/I_1 = \Psi \ot R/I_1$.
Since $\dim \CR = d$, applying (\ref{subtract})
 in the ring $\CR$ for the surjections $\Phi\ot
\CR$ and $\Psi\ot \CR$, we get a surjection
$\Gamma : P\ot \CR \surj I_2 \CR$ with 
$\Gamma \ot \CR/I_2
\CR = \Phi \ot \CR/I_2 \CR = \wt \phi \ot \CR/I_2 \CR.$

Again applying (\ref{subtract}) 
for the
surjections $\Gamma$ and $\wt \phi \ot \CR$, we get a surjection  
$\Gamma_1 : P\ot \CR \surj K\CR$ with $\Gamma_1 \ot \CR/K
\CR = \wt \phi \ot \CR/K \CR$.
Since $K+\CJ(A)=R$, applying (\ref{M2}), we get a
surjection $\Gamma_2 : P \surj K$ with $\Gamma_2 \ot R/K = \wt \phi
\ot R/K$.

We have two surjections $\wt \phi : P \surj I_2\cap K$ and $\Gamma_2 :
P \surj K$ with $\Gamma_2 \ot R/K = \wt \phi \ot R/K$.  Recall that
$K+\CJ(A)=R$. Following the proof of 
(\ref{addition-1}) Step 2, we get a surjection $\Delta : P \surj I_2$ with
$\Delta \ot R/I_2 = \wt \phi \ot R/I_2 = \Phi \ot R/I_2$.  This proves
the result.  $\hfill \gj$
\end{proof}

\begin{theorem}\label{M33}
Let $A$ be a ring of dimension $d$ and $R=A[T,T^{-1}]$. Let $n$ be an
integer such that $2n \geq d+3$.
Let $I$ be an ideal of $R$ of height 
$n$. Assume that $\hh \CJ(A) \geq d-n+2$. Let $P=P'\op R^2$ be a projective
$R$-module of rank $n$ and let $\phi : P \surj I/I^2$ be a
surjection. Assume that $\phi \ot \CR : P\ot \CR \surj I\CR/
I^2 \CR$ can be lifted to a
surjection $\Phi : P\ot \CR\surj I \CR$. Then $\phi$ can be
lifted to a surjection $\Delta : P \surj I$. 
\end{theorem}

\begin{proof}
Let $J=(I\cap A)\cap \CJ(A)$. Note
that $\hh J \geq d-n+2$.
Applying (\ref{moving}) with $f=1$, 
we get a lift $\Phi_1 \in \Hom(P,I)$ of
$\phi$ such that the ideal $\Phi_1(P)=I''$ satisfies the following
properties: 
$(1)~ I= I''+J^2$,
$(2)~ I'' = I\cap K$, where $\hh K\geq n$ and 
$(3)~ K+J^2=R$. 

If $\hh K > n$, then $K=R$ and $\Phi_1$ is a lift of $\phi$. 
Hence, we assume that $\hh K =n$. 
We have two surjections $\Phi : P\ot \CR \surj I\CR$ and $\Phi_1
: P \surj I\cap K$ with $\Phi \ot \CR/I\CR = \Phi_1 \ot
\CR/I\CR$. 
Applying (\ref{subtract}) in the
ring $\CR$ for the surjections $\Phi$ and $\Phi_1 \ot \CR$,
 we get a surjection $\Psi : P\ot \CR \surj K\CR$
such that $\Psi \ot \CR/K \CR = \Phi_1 \ot \CR/K \CR$.
Since $K+\CJ(A)=R$, applying
(\ref{M2}), we get a surjection $\Delta_1 : P \surj K$
which is a lift of $\Phi_1\ot R/K$.

We have two surjections $\Phi_1: P \surj  I\cap K$ and $\Delta_1 :
P \surj K$ with $\Phi_1 \ot R/K = \Delta_1 \ot R/K$.
Applying (\ref{subtract-1}), we get a surjection
$\Delta : P \surj I$ such that $\Delta \ot R/I = \Phi_1\ot R/I = \phi$.
This proves the result.
$\hfill \gj$
\end{proof}

As a consequence of the above result, we have the following:

\begin{corollary}\label{M3}
Let $A$ be a ring of dimension $n\geq 3$ with $\hh \CJ(A) \geq 2$ and
$R=A[T,T^{-1}]$. Let $I$ be an ideal of $R$ of height $n$. Let $\phi :
(R/I)^n \surj I/I^2$ be a surjection. Assume that $\phi\ot \CR$ can be
lifted to a surjection from $\CR^n$ to $I\CR$. Then $\phi$ can be
lifted to a surjection $\Phi : R^n \surj I$.
\end{corollary}

\section{Euler class group of $A[T,T^{-1}]$}

\notation We will denote the following hypothesis by {\bf (*)}:
Let $A$ be a ring containing $\BQ$
of dimension $n\geq 3$ with $\hh \CJ(A) \geq 2$ and
$R=A[T,T^{-1}]$.

Assume $(*).$ We proceed to define the $n^{th}$ {\it Euler class
group} of $R$. The results of this section are similar to
(\cite{Das}, Section 4), where it is proved for the ring $A[T]$ (without the
assumption $\hh \CJ(A) \geq 2$).

Let $I\subset R$ be an ideal of height $n$ such that $I/I^2$ is
generated by $n$ elements. Let $\ga$ and $\gb$ be two surjections from
$(R/I)^n$ to $I/I^2$. We say that $\ga$ and $\gb$ are {\it related} if there
exists $\sigma \in \SL_n(R/I)$ such that $\ga\sigma=\gb$. It is easy
to see that, this is an equivalence relation on the set of surjections
from $(R/I)^n$ to $I/I^2$. Let $[\ga]$ denote the equivalence class of
$\ga$. We call such an equivalence class $[\ga]$ a {\it local
orientation} of $I$.

If a surjection $\ga$ from $(R/I)^n$ to $I/I^2$ can be lifted to a
surjection $\gT : R^n \surj I$, then so can any $\gb$ equivalent to
$\ga$. For, let $\gb=\ga\sigma$ for some $\sigma\in \SL_n(R/I)$.  If
$I\CR=\CR$, then $\gb\ot \CR$ can be lifted to a surjection from
$\CR^n \surj I\CR$ and hence we can appeal to (\ref{M3}).
We assume that $I\CR$ is a proper ideal of $\CR$.
Since
$\dim \CR=n$, we have $\dim \CR/I\CR=0$. Hence,
$\SL_n(\CR/I\CR)=E_n(\CR/I\CR)$. Therefore, by 
(\ref{trans}), $\sigma \ot \CR$ can be lifted to an
element of $\SL_n(\CR)$. Thus
$\gb\ot \CR$ can
be lifted to a surjection from $\CR^n \surj I\CR$. By
(\ref{M3}), $\gb$ can be lifted to a surjection
from $R^n\surj I$. Therefore, from now on, we shall identify a
surjection $\ga$ with the equivalence class $[\ga]$ to which it
belongs. 

We call a local orientation $[\ga]$ of $I$ a {\it global orientation}
of $I$, if the surjection $\ga: (R/I)^n \surj I/I^2$ can be lifted to a
surjection $\gT : R^n \surj I$.

Let $G$ be the free abelian group on the set of pairs $(I,w_I)$, where
$I\subset R$ is an ideal of height $n$ having the property that $\Spec
(R/I)$ is connected, $I/I^2$ is generated by $n$ elements and $w_I:
R^n \surj I/I^2$ is a local orientation of $I$.

Let $I\subset R$ be an ideal of height $n$ such that $I/I^2$ is
generated by $n$ elements.. Then $I$ can be decomposed
as $I=I_1 \cap \ldots \cap I_r$, where $I_k$'s are pairwise comaximal
ideals of $R$ of height $n$ and $\Spec(R/I_k)$ is connected. From
(\cite{Das}, Lemma 4.4), it follows that such a decomposition is
unique. We say that $I_k$'s are the connected components of $I$. 
Let $w_I : (R/I)^n
\surj I/I^2$ be a surjection. Then $w_I$ induces surjections $w_{I_k}
: (R/I_k)^n \surj I_k/{I_k}^2$. By $(I,w_I)$, we denote the element $\sum
(I_k,w_{I_k})$ of $G$.

Let $H$ be the subgroup of $G$ generated by the set of pairs $(I,w_I)$,
where $I\subset R$ is an ideal of height $n$ and
$w_I$ is a global orientation of
$I$. We define the $n^{th}$ {\it Euler class group} of $R$, denoted by
$E^n(R)$, to be $G/H$. 
By abuse of notation, we will write $E(R)$ for
$E^n(R)$ throughout this paper.

Let
$P$ be a projective $R$-module of rank $n$ having trivial determinant.
Let $\gc : R\iso \wedge^n P$ be an isomorphism.
To the pair $(P,\gc)$, we associate an element $e(P,\gc)$ of $E(R)$ as
follows:

Let $\gl : P \surj I$ be a surjection, where $I\subset R$ is an ideal
of height $n$ (such a surjection exists by (\ref{EE})). 
Let ``bar'' denote reduction mod $I$.
We obtain an induced surjection $\ol \gl : P/IP \surj I/I^2$.
Since $P$ has trivial determinant and $\dim R/I \leq 1$, by
(\ref{serre}), $P/IP$ is a free $R/I$-module of rank $n$. We choose
an isomorphism $\ol \gamma :(R/I)^n \iso P/IP$ such that $\gw^n (\ol
\gamma) = \ol \gc$. Let $w_I$ be the surjection $\ol \gl\ol \gamma :
(R/I)^n \surj I/I^2$. Let $e(P,\gc)$ be the image of 
$(I,w_I)$ in $E(R)$. We say that $(I,w_I)$ is {\it obtained} from the
pair $(\gl,\gc)$.

\begin{lemma}
The assignment sending the pair $(P,\gc)$ to the element $e(P,\gc)$, as
described above, is well defined.
\end{lemma}

\begin{proof}
Let $\mu : P \surj I_1$ be another surjection, where $I_1\subset R$ is
an ideal of height $n$. Let $(I_1,w_{I_1})$ be obtained from the pair
$(\mu,\gc)$. Let $J=(I\cap I_1)\cap A$. Recall that $w_I : (R/I)^n
\surj I/I^2$ is a surjection.
By (\ref{moving}), $w_I$ can be lifted to
$\Phi:R^n \surj I\cap K$, where 
$\hh K=n$ and $K+J=R$. 

Since $K$ and $I$ are comaximal, $\Phi$ induces a local orientation
$w_K$ of $K$. Clearly, $(I,w_I)+(K,w_K)=0$ in $E(R)$.
Let $L=K \cap I_1$. Since $K+I_1=R$, $w_K$ and $w_{I_1}$ together induce a
local orientation $w_L$ of $L$. It is enough to show that $(L,w_L)=0$ in
$E(R)$ (Since
$(L,w_L)=(K,w_K)+(I_1,w_{I_1})$ in $E(R)$ and $(L,w_L)=0$ implies
$(I,w_I)=(I_1,w_{I_1})$ in $E(R)$).

Since $\dim \CR=n=$ rank $P$, $e(P\ot \CR,\gc \ot \CR)$ is
well defined in $E(\CR)$ (\cite{BR3}, Section 4). Hence, it follows
that $w_L\ot \CR$ is a global orientation of $L \CR$. Therefore,
by (\ref{M3}), $w_L$ is a global orientation of
$L$, i.e. $(L,w_L)=0$ in $E(R)$. This proves the lemma.  
$\hfill \gj$
\end{proof}

\notation We define the {\it Euler class} of $(P,\gc)$ to be $e(P,\gc)$.

\begin{theorem}\label{Zero}
Assume $(*)$.  Let
$I\subset R$ be an ideal of height $n$ such that $I/I^2$ is generated
by $n$ elements and let $w_I:R^n \surj I/I^2$ be a local orientation of
$I$. Suppose that the image of $(I,w_I)$ in 
$E(R)$ is zero. Then $w_I$ is a global orientation of $I$.
\end{theorem}

\begin{proof}
Since $(I,w_I)=0$ in $E(R)$, $(I\CR,w_I\ot \CR)=0$ in
$E(\CR)$. Therefore, by (\ref{bha}), $w_I\ot
\CR$ can be lifted to a surjection from $\CR^n \surj
I\CR$ (as $\dim \CR =n$). By (\ref{M3}), $w_I$ can be lifted to a
surjection from $R^n \surj I$ and hence is a global orientation of $I$.
$\hfill \gj$
\end{proof}

\begin{theorem}\label{manoj}
Assume $(*)$.
Let $P$ be a projective $R$-module 
of rank $n$ with trivial determinant
and let $I\subset R$ be an ideal of height $n$. Assume that, we are
given a surjection $\psi : P \surj I/I^2$. Assume further that, $\psi
\ot \CR$ can be lifted to a surjection $\Psi : P\ot \CR \surj I
\CR$. Then there exists a surjection $\wt \Psi : P \surj I$, which is a
lift of $\psi$. 
\end{theorem}

\begin{proof}
Let $J=I\cap \CJ(A)$. Then $\hh J \geq 2$.
By (\ref{moving}), $\psi$ can be lifted to $\Phi : P \surj I\cap I'$, where
$\hh I'=n$ and $I'+J^2=R$.

Fix $\gc : R \iso
 \wedge^n P$. Let  $\gl
 : (R/(I\cap I'))^n \iso P/(I\cap I')P$ such that $\wedge^n \gl = \gc
 \ot R/(I\cap I') $.
Then $e(P,\gc)=(I\cap I',w_{I\cap I'})$
 in $E(R)$, where $w_{I\cap I'} = (\Phi \ot R/(I\cap I'))\gl$.
Therefore, $e(P,\gc)=(I,w_I)+(I',w_{I'})$, where $w_I$ and
 $w_{I'}$ are local orientations of $I$ and $I'$ respectively induced
 from $w_{I\cap I'}$. 

Since $e(P\ot \CR,\gc \ot \CR) =(I\CR,w_I\ot \CR)$ (using $\Psi$),
$(I'\CR,w_{I'}\ot \CR)=0$ in $E(\CR)$, i.e. $w_{I'}\ot \CR$ can be
lifted to a surjection from $\CR^n$ to $I'\CR$.
 By (\ref{M3}), $w_{I'}$ can be lifted to $n$ set of generators of $I'$, say
$I'=(f_1,\ldots,f_n)$.
Since $I' +\CJ(A)=R$ and $\hh I'=n$, $\dim R/I'=0$.  Hence,
applying (\ref{lindel}, \ref{trans} and \ref{EE}); after performing
some elementary transformation on the generators of $I'$, we can
assume that

$(1)~ \hh (f_1,\ldots,f_{n-1})=n-1$,

$(2)~ \dim R/(f_1,\ldots,f_{n-1}) \leq 1$ and 

$(3)~ f_n=1$ modulo $J^2$.

Write $C=R[Y],~ K_1=(f_1,\ldots,f_{n-1},Y+f_n),~ K_2=IC$ and
$K_3=K_1\cap K_2$. \\

\noindent{\bf Claim :}
There exists a surjection $\Delta(Y) :
P[Y] \surj K_3$ such that $\Delta(0)=\Phi$.

First we show that the theorem follows from the claim. Specializing
$\Delta(Y)$ at $Y=1-f_n$, we obtain a surjection $\Delta_1 : P \surj
I$. Since $1-f_n \in J^2 \subset I^2$, $\Delta_1 = \Phi$ modulo $I^2$.
Therefore, $\Delta_1$ is a lift of $\psi$. This proves the result. \\

\noindent{\it Proof of the claim :}
$\gl$ induces an isomorphism $\gd : (R/I')^n \iso P/I'P$ such that
$\gw^n \gd =\gc \ot R/I'$. Also, $(\Phi \ot R/I')\gd = w_{I'}$.
Since $\dim C/K_1=\dim R/(f_1,\ldots,f_{n-1}) \leq 1$, and $P$ has
trivial determinant, by (\ref{serre}),
$P[Y]/K_1P[Y]$
is free of rank $n$. Choose an isomorphism $\Gamma(Y): (C/K_1)^n \iso
P[Y]/K_1P[Y]$ such that $\gw^n (\Gamma(Y)) = \gc \ot C/K_1$.   

Since $\gw^n \gd =\gc \ot R/I'$, $\Gamma(0)$ and $\gd$ differs 
by an element of $\SL_n(R/I')$. Since $\dim
R/I'=0$, $\SL_n(R/I')=E_n(R/I')$. Therefore, we can alter
$\Gamma(Y)$ by an element of $\SL_n(C/K_1)$ and assume that $\Gamma(0)
= \gd$.

Let $\Lambda(Y) : (C/K_1)^n \surj K_1/{K_1}^2$ be the surjection
induced by the set of generators $(f_1,\ldots,f_{n-1},Y+f_n)$ of
$K_1$. Thus, we get a surjection 
$$ \Delta(Y) =\Lambda(Y) \Gamma(Y)^{-1} : P[Y]/K_1P[Y] \surj K_1/{K_1}^2 .$$

Since $\Gamma(0)=\gd$, $\Phi \ot R/I'= w_{I'}\gd^{-1}$ and
$\Lambda(0) = w_{I'}$, we have $\Delta(0)=\Phi \ot R/I'$. 
By (\ref{Mandal-2}), we get a surjection $\wt \Delta : P[Y]
\surj K_3$ such that $\wt \Delta(0)=\Phi$. This proves the claim.
$\hfill \gj$
\end{proof}

\begin{lemma}\label{bm2}
Assume $(*)$. Let $P$
be a projective $R$-module of rank $n$ having trivial
determinant and $\gc :R \iso \gw^n P$. Let
$e(P,\gc) =(I,w_I)$ in $E(R)$, where $I\subset R$ is an ideal of height $n$.
Then there exists a surjection $\Delta :
P\surj I$ such that $(I,w_I)$ is obtained from $(\Delta,\gc)$. 
\end{lemma}

\begin{proof}
Since $\dim R/I \leq 1$ and $P$ has trivial determinant, by
(\ref{serre}), $P/IP$ is a free $R/I$-module of rank $n$. Choose $\gl :
(R/I)^n \iso P/IP$ such that $\gw^n \gl = \gc\ot R/I$.
Let $\gamma =w_I \gl^{-1} : P/IP \surj I/I^2$.

Since $e(P\ot \CR,\gc \ot \CR) = (I\CR,w_I\ot \CR)$ in
$E(\CR)$, by (\ref{bha}),
there exists a surjection $\Gamma : P\ot \CR \surj I\CR$ such that
$(I\CR,w_I\ot \CR)$ is obtained from the pair $(\Gamma, \gc \ot
\CR)$, i.e. $\Gamma$ is a lift of $\gamma \ot \CR$. Applying 
(\ref{manoj}), there exists a surjection $\Delta : P \surj I$ such
that $\Delta$ is a lift of $\gamma$. Since $(\Delta \ot
R/I)\gl=w_I$ and $\gw^n (\gl)=\gc\ot R/I$, $(I,w_I)$ is
obtained from the pair $(\Delta,\gc)$.
$\hfill \square$
\end{proof}

The following result is essentially (\ref{moving}).

\begin{lemma}\label{5.5}
Assume $(*)$.
Let $(I,w_I) \in
E(R)$. Then there exists an ideal $I_1\subset R$ of height $n$ and a
local orientation $w_{I_1}$ of $I_1$ such that $(I,w_I)+(I_1,w_{I_1})
= 0$ in $E(R)$. Further, $I_1$ can be chosen to be comaximal with any
ideal $K \subset R$ of height $\geq 2$.
\end{lemma}

\begin{corollary}\label{cor-1}
Assume $(*)$. Let $P$ be a projective $R$-module of rank $n$ with
trivial determinant and $\gc:R\iso \gw^n(P)$.
 Then $e(P,\gc)=0$ if and only if $P$ has a
unimodular element. In particular, if $P$ has a unimodular element,
then

$(1)$ $P$ maps onto any ideal of height $n$ generated by $n$ elements
(\ref{bm2}).

$(2)$ Let $\beta : P \surj I$ be a surjection, where $I$ is an ideal of
$R$ of height $n$.  Then $I$ is generated by $n$ elements.
\end{corollary}

\begin{proof}
Let $\ga : P \surj I$ be a surjection, where $I\subset R$ is an ideal
of height $n$. Let $e(P,\gc)=(I,w_I)$ in $E(R)$, where $(I,w_I)$ is
obtained from the pair $(\ga,\gc)$. 

Assume that $e(P,\gc)=0$ in $E(R)$. Then $(I,w_I)=0$ in
$E(R)$. By (\ref{5.5}), there exists an ideal $I'$ of height $n$ such
that $I'+\CJ(A)=R$ and
a local orientation $w_{I'}$ of $I'$ such that $(I,w_I)+(I',w_{I'})=0$
in $E(R)$. Since $(I,w_I)=0$,
$(I',w_{I'})=0$ in $E(R)$. Hence,
without loss of generality, we can assume that
$I+\CJ(A)R =R$.

By (\ref{Zero}), $I$ is generated by $n$ elements, say
$I=(f_1,\ldots,f_n)$. Since $I+\CJ(A) R =R$, $\dim R/I=0$. Hence,
applying (\ref{lindel}, \ref{trans}); after performing some
elementary transformations on the generators of $I$, we can assume
that $\dim R/(f_1,\ldots,f_{n-1})\leq 1$.

Let $C=R[Y]$ and $K=(f_1,\ldots,f_{n-1},Y+f_n)$ be an ideal of $C$. We
have two surjections $\ga : P \surj K(0) (=I)$ and $\phi : P[Y]/KP[Y]
\surj K/K^2$ such that $\phi(0) = \ga$ mod $K(0)^2$, where $\phi$ is
the composition of two maps, $\phi_1 : P[Y]/KP[Y] \iso (C/K)^n$ with
$\gw^n \phi_1 = \gc^{-1}\ot C/K$ and $\phi_2 : (C/K)^n \surj K/K^2$
defined by $(f_1,\ldots,f_{n-1},Y+f_n)$. Applying (\ref{Mandal-2})
with $I_1=K$ and $I_2=C$, we get a surjection $\Phi : P[Y] \surj
K$. Since $\Phi(1-f_n): P\surj R$, $P$ has a
unimodular element.

Conversely, we assume that $P$ has a unimodular element. Applying
(\ref{bha}), we have $(I\CR,w_I\ot
\CR)=0$ in $E(\CR)$. By ($\ref{M3}$), 
$(I,w_I)=0=e(P,\gc)$ in $E(R)$. This proves the result.
$\hfill \square$
\end{proof}

The following result is a direct consequence of (\ref{M3}). 

\begin{theorem}\label{manoj11}
Assume $(*)$.
Then the canonical map
$E(R)\ra E(\CR)$ is injective.
\end{theorem}

Assume $(*)$.  We have a canonical map $\Phi
: E(A)\ra E(R)$.  It is easy to see that $\Phi$ is injective. It is
natural to ask, when is $\Phi$ surjective? First, we prove an analogue
of (\cite{Bh-Manoj}, Theorem 4.13) for $A[T,T^{-1}]$.

\begin{theorem}\label{manoj-3}
Let $A$ be a regular domain of dimension $d$ essentially of finite
type over an infinite perfect field $k$ and $R=A[T,T^{-1}]$.  Let $n$
be an integer such that $2n \geq d+3$. Let $I \subset R$ be an ideal
of height $n$ and let $P$ be a projective $A$-module of rank
$n$. Assume that $I$ contains some $f\in A[T]$ such that either $f$ is
a monic polynomial or $f(0)=1$. Then any surjection 
$\phi : P\ot R \surj I/I^2$ can be lifted to a surjection $\Phi :
P\ot R \surj I$.
\end{theorem}

\begin{proof}
First we assume that $f(0)=1$.
Let $J=I \cap A[T]$.  Let $\psi : P\ot R \ra I$ be a lift of $\phi$.
Since $(P\ot R)^* = P^* \ot R$, there exists 
$\wt \psi \in P[T]^*$ such that $\psi = \wt \psi/T^r$ for some
positive integer $r$. It follows that $\wt \psi : P[T] \ra J$. 
Let $\Psi : P[T] \ra J/J^2$ be the map induced by $\wt \psi$.
Since $\Psi_T = \phi$ and $(J/J^2)_f=0$, we get that $\Psi$ is a
 surjection. Since $f\in I$, by (\cite{Bh-Manoj}, Lemma 3.5), $\Psi$
can be lifted to a surjection $\Delta : P[T] \surj J/J^2(f-1)$. Since
$f-1 \in (T)$, $\Delta$ induces a surjection $\wt \Delta : P[T] \surj
J/J^2T$.
Applying (\cite{Bh-Manoj}, Theorem 4.13), we get a
surjection $\Phi : P[T] \surj J$ which lifts $\wt \Delta$ and hence
$\Psi$. Now, $T^r(\Phi \ot R) : P\ot R \surj I$ is a
lift of $\phi$. This proves the result in the case $f(0)=1$.

Now, we assume that $f(T)$ is a monic polynomial. Let $J=I\cap A[X]$,
where $X=T^{-1}$. Then $J$ contains an element $g(X)=T^{-r}f(T)$,
where $r=$ deg $f$. Note that $g(0)=1$. Now, we are reduced to the
previous case.
$\hfill \square$     
\end{proof}

As a consequence of (\ref{manoj-3}), we have the following
result. 

\begin{theorem}\label{mandal-v}
Let $A$ be a regular domain of dimension $n\geq 3$ essentially of
finite type over an infinite perfect field $k$ with $\hh \CJ(A) \geq
2$.  Let $(I,w_I) \in E(A[T,T^{-1}])$. Assume that $I$ contains some
$f(T)\in A[T]$ such that either $f$ is a monic polynomial or
$f(0)=1$. Then $(I,w_I)=0$.
\end{theorem}

\begin{remark}
In \cite{mmm}, (\ref{manoj-3}) is proved for an arbitrary ring
under the assumption that $I$ contains a special monic
polynomial. Hence (\ref{mandal-v}) is valid for an arbitrary ring if
$I$ contains a special monic polynomial.
\end{remark}

Let $A$ be a ring of dimension $n$ containing an infinite field and
let $P$ be a projective $A[T]$-module of rank $n$. In
\cite{bhatwadekar}, it is proved that if $P_{f(T)}$ has a unimodular
element for some monic polynomial $f(T)\in A[T]$, then $P$ has a unimodular
element. We will prove the analogous result for $A[T,T^{-1}]$. 

\begin{theorem}
Assume $(*)$. Let $P$ be a projective $R$-module of rank $n$ with
trivial determinant. If $P_{f(T)}$ has a unimodular element for some
special monic polynomial $f(T)\in A[T]$, then $P$ has a unimodular
element.
\end{theorem}

\begin{proof}
Fix $\chi : R\iso \gw^n(P)$. Since $P_f$ has a unimodular element,
$e(P\ot \CR, \chi\ot \CR)=0$ in $E(\CR)$. By (\ref{manoj11}),
$e(P,\chi)=0$ in $E(R)$. Hence $P$ has a unimodular element, by
(\ref{cor-1}). $\hfill \square$
\end{proof}

\section{Weak Euler class group of $A[T,T^{-1}]$}

Results in this section are similar to (\cite{Das}, Section 5).
Assume $(*)$. We define the $n^{th}$ {\it weak Euler class group}
${E_0}^n(R)$ of $R$ in the following way :

Let $G$ be the free abelian group on $(I)$, where $I\subset R$ is an
ideal of height $n$ with the property that $I/I^2$ is generated by $n$
elements and $\Spec(R/I)$ is connected.  Let $I\subset R$ be an ideal
of height $n$ such that $I/I^2$ is generated by $n$ elements. Then $I$
can be decomposed as $I=I_1\cap \ldots \cap I_r$, where $I_i$'s are
pairwise comaximal ideals of height $n$ and $\Spec (R/I_i)$ is
connected for each $i$. In the previous section, we have seen that
such a decomposition of $I$ is unique.  By $(I)$, we denote the
element $\sum_i (I_i)$ of $G$.

Let $H$ be the subgroup of $G$ generated by elements of the type
$(I)$, where $I\subset R$ is an ideal of height $n$ such that $I$ is
generated by $n$ elements. 

We define ${E_0}^n(R)=G/H$.
By abuse of notation, we will write $E_0(R)$ for ${E_0}^n(R)$ in what
follows. Note that, there is a canonical surjective homomorphism from
$E(R)$ to $E_0(R)$ obtained by forgetting the orientations.

\begin{remark} 
Assume $(*)$.  Let $I\subset R$ be an ideal of height $n$ and let $w_I
: (R/I)^n \surj I/I^2$ be a local orientation of $I$. Let $\gt \in
\GL_n(R/I)$ be such that det $\gt = \ol f$. Then $w_I \gt$ is another
orientation of $I$, which we denote by $\ol f w_I$. On the other hand,
if $w_I$ and $\wt w_I$ are two local orientations of $I$, then by
(\cite{BR3}, Lemma 2.2), it is easy to see that $\wt w_I = \ol f w_I$
for some unit $\ol f\in R/I$.
\end{remark}
\medskip

The proof of the following lemma is contained in (\cite{BR3},
2.7, 2.8 and 5.1) and hence, we omit the proof.

\begin{lemma}\label{5.2}
Assume $(*)$. Let $P$ be a projective
$R$-module of rank $n$ having trivial determinant and $\gc : R\iso
\gw^n P$. Let $\ga : P \surj I$ be a surjection,
where $I\subset R$ is an ideal of height $n$. Let $(I,w_I)$ be
obtained from $(\ga,\gc)$. Let $f\in R$ be a unit mod $I$. Then
there exists a projective $R$-module $P_1$ of rank $n$ 
such that $[P]=[P_1]$ in $K_0(R)$, 
$\gc_1 : R\iso \gw^n P_1$ and a surjection
$\gb : P_1 \surj I$ such that $(I, \ol {f^{n-1}} w_I)$
is obtained from $(\gb,\gc_1)$.
\end{lemma}

The following lemma can be proved using (\cite{BR3}, Lemma 5.3, 5.4)
and (\ref{M3}).

\begin{lemma}\label{5.3}
Assume $(*)$. 
Let $(I,w_I)\in E(R)$.  Let $\ol f
\in R/I$ be a unit. Then $(I,w_I)=(I,\ol {f^2}w_I)$ in $E(R)$.
\end{lemma}

Adapting the proof of (\cite{Bh-Raja-2}, Lemma 3.7) and using
(\ref{EE}) in place of Swan's Bertini
theorem, the proof of the following lemma follows.

\begin{lemma}\label{5.6}
Assume $(*)$ with $n$ even. Let $P$ be a stably
free $R$-module of rank $n$ and $\gc : R\iso \gw^n P$. Suppose that
$e(P,\gc)=(I,w_I)$ in $E(R)$. Then $(I,w_I)=(I_1,w_{I_1})$ in $E(R)$
for some ideal $I_1 \subset R$ of height $n$ generated by $n$
elements. Moreover, $I_1$ can be chosen to be comaximal with any ideal
of $R$ of height $\geq 2$.
\end{lemma}

The following result can be proved by adapting the proofs
of (\cite{Bh-Raja-2}, 3.8, 3.9, 3.10, 3.11).

\begin{proposition}\label{5.7}
Assume $(*)$ with $n$ even. Then we have the
followings: 

$(1)$ Let $I_1,
I_2 \subset R$ be two comaximal ideals of height $n$ and
$I_3=I_1 \cap I_2$. If any two of $I_1,I_2$ and $I_3$ are surjective images of
stably free  $R$-modules of rank $n$, then so is the third.

$(2)$ Let $(I,w_I) \in E(R)$. Then $(I)=0$ in $E_0(R)$ if and only if
$I$ is a surjective image of a stably free projective $R$-module of
rank $n$.

$(3)$ Let $P$ be a projective $R$-module of rank $n$ with trivial
determinant.  Then $e(P)=0$ in $E_0(R)$ if and only if $[P]=[Q\op R]$
in $K_0(R)$ for some projective $R$-module $Q$ of rank $n-1$.

$(4)$ Let $P$ be a
projective $R$-module of rank $n$ with trivial determinant. Suppose
that $e(P)=(I)$ in $E_0(R)$, where $I\subset R$ is an ideal of height
$n$. Then there exists a projective $R$-module $Q$ of rank $n$ such
that $[Q]=[P]$ in $K_0(R)$ and $I$ is a surjective image of $Q$.
\end{proposition}

The proof of the following result is same as of (\cite{BR3},
Proposition 6.5) using above results.

\begin{theorem}\label{5.10}
Assume $(*)$ with $n$ even. Let $(I,\wt
{w_I})\in E(R)$ belongs to the kernel of the canonical homomorphism
$E(R)\surj E_0(R)$. Then there exists a stably free $R$-module $P_1$
of rank $n$ and $\gc_1 : R\iso \gw^n P_1$ such that
$e(P_1,\gc_1)=(I,\wt {w_I})$ in $E(R)$.
\end{theorem}


\section{The case of dimension two}

In this section, we briefly outline the results similar to those in the
previous sections in the case when dimension of the base ring is
two. The results of this section are similar to (\cite{Das}, Section
6), where it is proved for $A[T]$.

We begin by stating the following result of Mandal (\cite{Mandal}).

\begin{lemma}\label{mandal91}
Let $A$ be a ring and $R=A[T,T^{-1}]$. Let $P$ be a projective
$R$-module. Let $f\in R$ be a special monic polynomial. If $P_f$ is
free, then $P$ is free.
\end{lemma}

The proof of the following result is similar  to
(\cite{Das}, Theorem 7.1). 

\begin{theorem}\label{dim2}
Let $A$ be a ring of dimension $2$ and $R=A[T,T^{-1}]$. Let $I\subset
R$ be an ideal of height $2$ such that $I=(f_1,f_2)+I^2$. Suppose that
there exists $F_1,F_2 \in I\CR$ such that $I\CR=(F_1,F_2)$ and
$F_i -f_i \in I^2\CR$ for $i=1,2$.  Then there exists $h_1,h_2\in
I$ and $\gt \in \SL_2(R/I)$ such that $I=(h_1,h_2)$ and $(\ol f_1, \ol
f_2)\gt = (\ol h_1,\ol h_2)$, where ``bar'' denotes reduction modulo
$I$.
\end{theorem}

\begin{proof}
Since a unimodular row of length two is always completable to a matrix
of determinant $1$, it follows (using patching
argument) that there is a projective $R$-module $P$ of rank $2$ 
with trivial
determinant mapping onto $I$. Let $\ga : P \surj I$ be the
surjection. Fix $\gc:R\iso \gw^2 P$. Since $\dim R/I
\leq 1$, by (\ref{serre}), $P/IP$ is free of rank $2$. Hence $\ga$ 
and $\gc$ induces a set of generators of $I/I^2$, say
$I=(g_1, g_2)+I^2$.

It is easy to see that there exists a matrix $\ol \sigma \in
\GL_2(R/I)$ with determinant $\ol f$ such that $(\ol f_1,\ol f_2) =
(\ol g_1,\ol g_2)\ol \sigma$. Now, following (\cite{BR3}, Lemma
2.7 and Lemma 2.8), there exists a projective $R$-module
$P_1$ of rank $2$ having trivial determinant, $\gc_1: R\iso
\gw^2 P_1$ and a surjection $\gb: P_1 \surj I$ such that if the
set of generators of $I/I^2$ induced by $\gb$ and $\gc_1$ is $\ol
h_1,\ol h_2$, then $(\ol h_1, \ol h_2)=(\ol g_1,\ol g_2)\ol \gd$,
where $\ol \gd \in \GL_2(R/I)$ has determinant $\ol f$.
Therefore, the two set of generators, $(\ol f_1, \ol
f_2)$ and $(\ol h_1,\ol h_2)$ of $I/I^2$ are connected by a matrix in
$\SL_2(R/I)$.

From the above discussion, it is easy to see  that $e(P_1\ot \CR,\gc_1\ot
\CR) = (I \CR,w_I \ot \CR)$ in $E(\CR)$, where $w_I :
(R/I)^2 \surj I/I^2$ is the surjection corresponding to the generators
$(\ol f_1,\ol f_2)$. Therefore, from the given condition of the
theorem, it follows that $(I \CR,w_I \ot \CR) = 0$ in
$E(\CR)$. Hence, we have $e(P_1\ot \CR,\gc_1\ot
\CR) = 0$ in $E(\CR)$. Since $\dim \CR=2$, by 
(\ref{bha}), 
$P_1 \ot \CR$ has a unimodular element and
hence is free (as rank $P_1=2$ and determinant of $P_1$ is trivial).
 Therefore, by (\ref{mandal91}), $P_1$ is a free $R$-module. 

Assume that the surjection $\gb$ is given by $h_1,h_2$. Then
$I = (h_1,h_2)$ and $(\ol f_1,\ol f_2) \gt = (\ol h_1,\ol h_2)$, for
some $\gt \in \SL_2(R/I)$. This proves the result.
$\hfill \square$
\end{proof}

As applications of the above theorem, we prove the following results.

\begin{corollary}
{\rm (Addition Principle)}
Let $A$ be a ring of dimension $2$ and $R=A[T,T^{-1}]$. Let $I_1,I_2\subset
R$ be two comaximal ideals of height $2$. Suppose that
$I_1=(f_1,f_2)$ and $I_2=(g_1,g_2)$.
Then there exists $h_1,h_2\in I_1\cap I_2$ and $\gt_i \in
\SL_2(R/I_i)$, $i=1,2$,
such that $I_1\cap I_2=(h_1,h_2)$ and $(f_1,f_2)\ot R/I_1 = ((h_1,
h_2)\ot R/I_1) \gt_1$ and $(g_1,g_2)\ot R/I_2 = ((h_1,h_2)\ot
R/I_2)\gt_2$. 
\end{corollary}

\begin{proof}
Write $I$ for $I_1\cap I_2$.
The generators of $I_1$ and $I_2$ induce a set of generators of
$I/I^2$, say $I=(H_1,H_2)+I^2$. Since $\dim \CR=2$, 
applying 
(\ref{addition}) in the ring $\CR$, we get $I\CR = (F_1,F_2)$
with $F_i - f_i\in {I_1}^2 \CR$ and $F_i - g_i \in {I_2}^2
\CR$. Hence, it is easy to see that $F_i - H_i \in I^2 \CR$, for
$i=1,2$.

Applying (\ref{dim2}), there exists $h_1,h_2\in I$
and $\gt \in \SL_2(R/I)$ such that $I=(h_1,h_2)$ and $(H_1,H_2)\ot R/I
= ((h_1,h_2)\ot R/I)\gt$. Let $\gt_i = \gt \ot R/I_i$. Then
$\gt_i\in \SL_2(R/I_i)$, $i=1,2$ and we have $(f_1,f_2)\ot R/I_1 =
((h_1, h_2)\ot R/I_1) \gt_1$ and $(g_1,g_2)\ot R/I_2 = ((h_1,h_2)\ot
R/I_2)\gt_2$.  $\hfill \square$
\end{proof}

\begin{corollary}
{\rm (Subtraction Principle)}
Let $A$ be a ring of dimension $2$ and $R=A[T,T^{-1}]$. Let $I_1,I_2\subset
R$ be two comaximal ideals of height $2$. Suppose that
$I_1=(f_1,f_2)$ and $I_1\cap I_2=(h_1,h_2)$ such that $f_i-h_i \in
{I_1}^2$, for $i=1,2$. Then there exists $g_1,g_2\in I_2$ and $\gt \in
\SL_2(R/I_2)$ such that $I_2=(g_1,g_2)$ and $(g_1,g_2)\ot R/I_2 =
((h_1,h_2)\ot R/I_2)\gt$.
\end{corollary}

\begin{proof}
We have $I_2=(h_1,h_2)+{I_2}^2$. Since $\dim \CR=2$, 
applying 
(\ref{subtract}) in the ring $\CR$, we get that $I_2 \CR =
(G_1,G_2)$ with $G_i-h_i \in {I_2}^2 \CR$. 
Now, applying (\ref{dim2}), we get the result.
$\hfill \square$
\end{proof}

\begin{remark}
Let $A$ be a ring of dimension $2$ and $R=A[T,T^{-1}]$. 
We can define the {\it
Euler class group} and the {\it weak Euler class group} of
$R$ in exactly the same way as we did in the previous
sections. The only difference is that, for an ideal $I$ of $R$ of
height $2$, a local orientation $[\ga]$ will be called a global
orientation if there is a surjection $\gt : R^2 \surj I$ and some
$\sigma \in \SL_2(R/I)$ such that $\ga \sigma = \gt \ot R/I$. For a
rank $2$ projective $R$-module $P$ having trivial determinant, the
Euler class of $P$ is defined as in the previous section.
\end{remark}

The following result can be proved using (\ref{dim2}, \ref{bha}) ($(i)$
follows from (\ref{Zero}), $(ii)$'s
proof is similar to (\cite{Das}, Theorem 7.6) using
(\ref{88}) and $(iii,iv)$ follows from (\ref{mandal91})).

\begin{theorem}
Let $A$ be a ring of dimension $2$ and $R=A[T,T^{-1}]$. Let $I\subset
R$ be an ideal of height $2$ such that $I/I^2$ is generated by $2$
elements. Let $w_I : (R/I)^2 \surj I/I^2$ be a local orientation of
$I$. Let $P$ be a projective $R$-module of rank $2$ with trivial
determinant and $\gc:R\iso \gw^2 P$. We have the following results:

$(i)$ Suppose that the image of $(I,w_I)$ is zero in
$E(R)$. Then $I$ is generated by $2$ elements and
$w_I$ is a global orientation of $I$.

$(ii)$ Suppose that $e(P,\gc)=(I,w_I)$ in $E(R)$. Then there exists a
surjection $\ga : P \surj I$ such that $(I,w_I)$ is obtained from
$(\ga,\gc)$.

$(iii)$ $e(P,\gc)=0$ in $E(R)$ if and only if $P$ has a unimodular
element  and hence $P$ is free.

$(iv)$ The canonical map $E(R)\ra E(\CR)$ is injective.
\end{theorem}

\begin{remark}
Let $A$ be a ring of dimension $2$ and  $R=A[T,T^{-1}]$. Let
$I\subset R$ be an ideal of height $2$ such that $I/I^2$ is generated
by $2$ elements and let $w_I$ be a local orientation of $I$. It is
easy to see, as in (\ref{dim2}), that there exists a projective
$R$-module $P$ of rank $2$ together with an isomorphism $\gc : R\iso
\gw^2 P$ and a surjection $\ga : P \surj I$ such that $(I,w_I)$ is
obtained from the pair $(\ga,\gc)$

The theory of weak Euler class group described in the last section
also follows in a like manner in the two dimensional case. 
\end{remark}


\section{Relations Between $E(R)$ and $\wt K_0Sp(R)$}

In this section, we prove results similar to (\cite{BR3},
Section 7).

Let $A$ be a ring of dimension $2$ and $R=A[T,T^{-1}]$. Let $\wt K_0
Sp(R)$ be the set of isometry classes of $(P,s)$, where $P$ is a
projective $R$-module of rank $2$ with trivial determinant
and $s:P\times P \ra R$ a
non-degenerate alternating bilinear form. We note that there is (up-to
isometry) a unique non-degenerate alternating bilinear form on $R^2$,
which we denote by $h$, namely $h((a,b),(c,d))=ad-bc$. We write $H(R)$ for
$(R^2,h)$.

We define a binary operation $*$ on $\wt K_0 Sp(R)$ as follows. Let
$(P_1,s_1)$ and $(P_2,s_2)$ be two elements of $\wt K_0 Sp(R)$.  Since
$\dim A=2$, $R=A[T,T^{-1}]$ and $P_1\op P_2$ has rank $4$, hence by
(\ref{serre}), $P_1\op P_2$ has a unimodular element, say $p$. Then
there exists $q\in P_1\op P_2$ such that if $s=s_1\perp s_2$, then
$s(p,q)=1$. Let $P_3=\{\wt p \in P_1\op P_2 \,|\, s(p,\wt p)=0=s(q,\wt
p)\}$. Then the restriction $s_3 : P_3 \times P_3 \ra R$ of $s$ to
$P_3$ is non-degenerate (i.e. $(P_3,s_3)$ is symplectic) and $P_1\op
P_2=(Rp\op Rq)\op P_3$.  Hence $(P_1,s_1)\perp (P_2,s_2)$ is isometric
to $(P_3,s_3)\perp (R^2,h)$. We define $(P_1,s_1) * (P_2,s_2)=(P_3,s_3)$.
By (\ref{88}), $(P_3,s_3)$ is determined
uniquely up-to isometry.
Hence $*$ is well defined operation and for
every symplectic $R$-module $(P,s)$ of rank $2$,
$(P,s)*(R^2,h)=(P,s)$. Hence $\wt K_0Sp(R)$ is a commutative semigroup
under $*$ with the isometry class of $(R^2,h)$ as the identity
element. We will briefly indicate that infact $\wt K_0
Sp(R)$ is an abelian group under $*$.

For a projective $R$-module $P$ of rank $2$ with trivial determinant,
the alternating bilinear form $s_P$ on $P\op P^*$ defined by
$$s_P((p,f),(q,g))=g(p)-f(q), \;p,q\in P,\;f,g\in P^*$$ is
non-degenerate. We write $H(P)$ for the symplectic module $(P\op
P^*,s_P)$.  If $(P,s)$ is a symplectic $R$-module of rank $2$, then
$(P,s)\perp (P,-s) \iso H(P)$ (\cite{sw4}, Lemma A.3). By
(\cite{mandal1}, Theorem 2.1), every projective $R$-modules of rank
$\geq 3$ has a unimodular element. Hence, by (\ref{B}), there exists a
projective $R$-module $P_1$ of rank $2$ such that $P\op P_1 \iso R^4$.
Therefore 
$$H(P_1)\perp (P,-s)\perp (P,s) \iso H(P_1\op P)\iso H(R^4)
\iso H(R^2)\perp H(R)\perp H(R).$$  
Since the symplectic module
$H(P_1)\perp (P,-s)$ has rank $6$, $H(P_1)\perp (P,-s) \iso
H(R^2)\perp (\wt P,\wt s)$ for some symplectic $R$-module $(\wt P,\wt
s)$ of rank $2$. By Bass result \cite{bass2}, 
$$(\wt P,\wt s)\perp
(P,s) \iso H(R)\perp H(R)$$ 
and therefore $(\wt P,\wt
s)*(P,s)=H(R)$. Thus, $\wt K_0Sp(R)$ is an abelian group under $*$.

Let $P$ be a projective $R$-module of rank $2$ with trivial
determinant. Then having a
non-degenerate alternating bilinear form $s$ on $P$ is equivalent to giving an
isomorphism $\gl : \gw^2 P\iso A$. Thus, we can identify the pair
$(P,s)$ with $(P,\gc)$, where $\gc$ is the generator of $\gw^2 P$
given by $\gl^{-1}(1)$. It is easy to see that the isometry classes of
$(P,s)$ coincides with the isomorphism classes of $(P,\gc)$.

We will begin with the following result, the proof of which 
is same as of (\cite{BR3},
Theorem 7.2).

\begin{theorem}
Let $A$ be a ring of dimension $2$ and $R=A[T,T^{-1}]$. Then the map
from $\wt K_0 Sp(R)$ to $E(R)$ sending $(P,\gc)$ to $e(P,\gc)$ is an
isomorphism. 
\end{theorem}

Let $A$ be a ring of dimension $2$ and $R=A[T,T^{-1}]$. Let $G$ be the
set of isometry classes of non-degenerate alternating bilinear forms
on $R^4$. Let $H(R^4) =(R^2,h)\perp (R^2,h)$. As before, we can define
the group structure on $G$ as follows: We set $(R^4,s_1) * (R^4,s_2) =
(R^4,s_3)$, where $s_3$ is the unique (up-to isometry) alternating
bilinear form on $R^4$ satisfying the property that $(R^4,s_1) \perp
(R^4,s_2)$ is isometric to $(R^4,s_3) \perp H(R^2)$. Then $G$ is a
group with $H(R^2)$ as the identity element. Let $s$ be a
non-degenerate alternating bilinear form on $R^4$. Since $\dim
A=2$ and $R=A[T,T^{-1}]$, by(\ref{B}), we get $(R^4,s)\iso (P,s') \perp
(R^2,h)$. The assignment sending $(R^4,s)$ to $(P,s')$ gives rise to
an injective homomorphism from $G$ to $\wt K_0 Sp(R)$.

In view of the above theorem, we have the following result, the proof
of which is same as (\cite{BR3}, Theorem 7.3).

\begin{theorem}
Let $A$ be a ring of dimension $2$ and $R=A[T,T^{-1}]$. Then we have
the following exact sequence
$$0 \ra G \ra \wt K_0 Sp(R) (\iso E(R)) \ra E_0(R) \ra 0.$$
\end{theorem}

\begin{corollary}\label{manoj14}
Assume $(*)$. Let
$(I,w_I)$ be an element of $E(R)$ such that its image in $E_0(R)$
(which is independent of $w_I$) is zero. Then the element
$(I,w_I)+(I,-w_I)=0$ in $E(R)$.
\end{corollary}

\begin{proof}
Let $(I,w_I)+(I,-w_I)=(J,w_J)$ in $E(R)$. 
Since $\dim \CR=n$, applying (\cite{BR3}, Corollary 7.9) in
the ring $\CR$, we get that $(J\ot \CR,w_J\ot \CR)=0$ in
$E(\CR)$. By (\ref{manoj11}),
$(J,w_J)=0$ in $E(R)$. This proves the result.
$\hfill \square$
\end{proof}

As an application of (\ref{manoj14}), following the proof of
(\cite{BR3}, Corollary 7.10), we have the following result.

\begin{corollary}
Assume $(*)$ with $n$ odd. Let $P$ be a
projective $R$-module of rank $n$ having trivial determinant. Assume
that the kernel of the canonical surjection $E(R) \surj E_0(R)$ has no
non-trivial 2-torsion. If $e(P)=0$ in $E_0(R)$, then $P$ has
a unimodular element.
\end{corollary}



Following the proof of (\cite{BR3}, Theorem 7.13) gives the following
result.

\begin{theorem}
Assume $(*)$ with $n$ odd. Let $P$ be a
projective $R$-module of rank $n$ having trivial determinant. Suppose
that there exists a projective $R$-module $Q$ of rank $n-1$ such that
$[P] =[Q\op R]$ in $K_0(R)$. Then $P$ has a unimodular element.
\end{theorem}


\section{Appendix}

We will freely use results and notations from \cite{Bh-2}. 
Let $(P,\langle,\rangle)$ be an $A$-module with an alternating
bilinear form $\langle,\rangle$ ($P$ need not be projective and
$\langle,\rangle$ need not be non-degenerate). 
Let $E(A^2\perp P,\langle,\rangle)$ denote the subgroup of $\Aut(A^2\perp P,
\langle,\rangle)$ generated by $\gt_{(c,q)}$ and $\sigma_{(d,q)}$ for $c,d\in
A$ and $q\in P$, where $\gt_{(c,q)}$  and $\sigma_{(d,q)}$ are defined as
$$\gt_{(c,q)}(a,b,p) = (a,b+ca+\langle p,q\rangle, p+aq),$$ 
$$ \sigma_{(d,q)}  (a,b,p) =
(a+bd+\langle q,p\rangle,b,p+bq)$$ for $(a,b,p)\in A^2\op P$.

\begin{remark}
It is easy to see that
(\cite{Bh-2}, Lemma 4.3, 4.5 4.7) holds for $(P,\langle,\rangle)$
replacing $ESp(A^2\perp P,\langle,\rangle)$ with $E(A^2\perp
P,\langle,\rangle)$ with
further assumption in (4.5) that $sP\subset F$. 
\end{remark}

The following result is a symplectic analogue of (\ref{lindel}) and is
a generalization of \cite{bass2} and
(\cite{Bh-2}, Theorem 4.8), where it is proved for $r=r'=0$ and $r=0$
respectively.  Our proof closely follows \cite{Bh-2}.

\begin{theorem}\label{88}
Let $B$ be a ring of dimension $d$ and $A=B[Y_1,\ldots,Y_{r'},X_1^{\pm
1},\ldots, X_r^{\pm 1}]$. Let $(P,\langle,\rangle)$ be a symplectic $A$-module
of rank $2n > 0.$ If $2n \geq d$, then $ESp(A^2 \perp P, \langle,\rangle)$ acts
transitively on $\Um(A^2 \oplus P).$
\end{theorem}

\begin{proof}
Let $(g_1, g_2, p) \in \Um(A^2 \oplus P)$. We want to show that there
exists $\Gamma \in ESp(A^2 \perp P, \langle,\rangle)$ such that $\Gamma
(g_1, g_2, p) = (1,0,0).$ We prove the result by induction on $r$.

If $r=0$, then the result follows from (\cite{Bh-2}, Theorem 4.8). Hence, we
assume that the result is proved for $r-1$ and $r\geq 1$. For the sake of
simplicity, we write $R=B[Y_1,\ldots,Y_{r'}, X_1^{\pm 1},\ldots,X_{r-1}^{\pm
1}]$ and $X_r=X$.

Without loss of generality, we can assume that $B$ is reduced. Let $S$
be a set of non-zero-divisors of $B$. Then $B_S$ is a finite direct
product of fields and therefore, by \cite{Suslin-4, Sw}, 
every projective
$A_S$-module is free. Hence, we can find a basis $\wt {p_1}, \ldots,
\wt {p_n}, \wt {q_1}, \ldots, \wt {q_n}$ of $P_S$ such that $\langle
\wt {p_i},\wt {p_j}\rangle = 0 = \langle \wt {q_i},\wt {q_j}\rangle$ for
$1 \leq i,j \leq n$, $\langle \wt {p_i},\wt {q_i}\rangle = 1$ for $1\leq i\leq
n$ and
$\langle \wt {p_i},\wt {q_j}\rangle = 0$ for $1\leq i,j \leq n$, $i \not =
j$. 

We can choose some $t\in S$ such that 
$\wt {p_i} = e_i/t$, $\wt {q_i} = f_i/t$ for some  $e_i, f_i \in
P$ for $1 \leq i \leq n $. Let $s = t^2$ and  $F = \sum_{i=1}^n Ae_i
+ \sum_{i=1}^n Af_i$. Then, by (\cite{Bh-2}, Lemma 4.2), 
$F$ is a free $A$-submodule of
$P$ of rank $2n$ and $sP \subset F$.

Let $F_1 = \sum_{i=1}^{n} R[X]e_i + \sum_{i=1}^{n} R[X]f_i$.
Let $P$ be generated by $\mu_1, \cdots,\mu_l$ as an $A$-module such that 
$(1)$ the set $\mu_1,\cdots,\mu_l$ contains $e_1,\cdots,e_n,f_1,\cdots,f_n$,
$(2)$ $s\mu_i \in F_1$ for $1\leq i\leq l$
and $(3)$ $\langle \mu_i,\mu_j\rangle \in R[X]$ for $1 \leq i,j \leq l$. Let
$M=\sum_{i=1}^l R[X] \mu_i$. Then $MA=P$ and $sM \subset F_1$.

Since $s \in B$ is a non-zero-divisor, $B_1 = B[X^{\pm 1}]/(s(X-1))$ is a
ring of dimension $d$ and $\ol A = A/(s(X-1))=B_1[Y_1,\ldots,Y_{r'}, X_1^{\pm
    1}, \ldots,X_{r-1}^{\pm 1}]$. Moreover, since rank $P \geq d$, by
(\cite{L}, Theorem 1.19), 
the map $\Um (A^2 \oplus P) \ra \Um (\ol
{A^2} \oplus P/s(X-1)P) $ is surjective. Therefore, by (\cite{Bh-2}, Lemma
4.1) and induction hypothesis,
there exists $\Psi \in ESp(A^2 \perp P, \langle ,\rangle)$ such that
$\Psi (g_1, g_2, p)= (1,0,0)$ modulo $s(X-1)A$.

Replacing $(g_1,g_2,p)$ with $\Psi (g_1, g_2, p)$, we may assume that
 $(g_1,g_2,p)=(1,0,0)$ modulo $s(X-1)A$.
By (\ref{EE}), there exist $h \in A$ and $p_1 \in P$
 such that $\hh(Ag_3 + I) \geq  d+1$,
 where $g_3 = g_1 + h g_2$, $p_2 = p +
 g_2 p_1$ and $I = p_2(P^*)= \langle P,p_2 \rangle$. Put
 $ \alpha = g_3 + \langle p_1, p\rangle \in A$. Then
 $$\sigma_{(h,p_1)}(g_1,g_2,p) = (g_1+g_2h+\langle p_1,p\rangle, g_2,p+g_2p_1)
=  (\alpha,g_2,p_2).$$
Since $(g_3,g_2,p)=(1,0,0)$ modulo $s(X-1)A$,
$\alpha =1$ modulo $s(X-1)A$.  Moreover, since
$\langle p_1,p_2\rangle = \langle p_1, p\rangle  \in I$, 
$(g_3, I)A = (\alpha, I)A
= (\alpha, s(X-1)I)A$.
Now, since $ (\alpha, s(X-1)I)A$ is an ideal of $A$ of
height $> d = \dim B$, by Mandal's theorem \cite{mandal-s}, $ (\alpha,
s(X-1)I)A$ contains a special monic polynomial, say $\gamma$, in the variable
$X$. We write $\gamma=\gamma(X)\in R[X]$.

Let $\beta(X) = g_2 + \gamma(X) \gamma_1$ for some suitable
$\gamma_1 \in A$ such that $\beta(X) \in R[X]$ and is a special monic
polynomial. Let $\gamma(X)\gamma_1 = \mu \alpha + \nu$ for some $\mu \in A$ and
$\nu \in s(X-1)I$. Since $I = \langle P,p_2\rangle$, there
exists $p_3 \in s(X-1)P$ such that $\nu = \langle -p_3, p_2\rangle
= \langle p_2,p_3\rangle$. Put $p_4 = p_2 + \alpha p_3.$ Then
$$\theta_{(\mu, p_3)}(\alpha,g_2, p_2)  =(\alpha,g_2+\mu\alpha+\langle
p_2,p_3\rangle,p_2 +\alpha p_3)=(\alpha,
\beta(X), p_4).$$ 
Note that, $(\alpha,p_4)=(1,0)$ modulo $s(X-1)A$.
and $\beta(X)$ is special monic polynomial.

Since $sP\subset F$, let $p_4 = (X-1) (\sum_{i=1}^n h_i e_i +
\sum_{j=1}^n k_j f_j)$ for some $h_i,k_j\in A$. Let $h_1 = -\gl X^{-r_0}
+ \wt h_1$, where $\wt h_1 \in A$ has $X^{-1}$ degree $\leq r_0-1$ and
$\gl \in R$. Let $a_0= (X-1)X^{-r_0}\gl$.
Then $$\sigma_{(0,a_0e_1)}(\alpha, \beta(X), p_4) =
(\alpha + a_0 \langle e_1,p_4\rangle , \beta(X), p_4 +
\beta(X) a_0e_1).$$  Note that, if $p_4 +
a_0\beta(X) e_1 = (X-1)(e_1 h_{11} + \sum_{i=2}^n h_i e_i +
\sum_{j=1}^n k_j f_j)$, then degree of $X^{-1}$ in $h_{11}\in A$ is
$\leq r_0-1$. Also note that $\alpha + a_0 \langle
e_1,p_4\rangle =1$ modulo $s(X-1)A$.  Hence, by induction on the $X^{-1}$
degree, applying such
symplectic transvections, say $\Psi_1 \in ESp(A^2 \perp P, \langle
,\rangle)$, we can assume that if $\Psi_1(\ga,\gb(X),p_4)=(\ga_1,\gb(X),
p_5)$, then
$p_5 \in (X-1)F_1$. Now, we write $p_5$ as $
p_5(X)$. We still have $\ga_1=1$ mod $s(X-1)A$. Write $\Gamma_1 =\Psi_1
\theta_{(\mu, p_3)} \sigma_{(h,p_1)}$. Then $\Gamma_1
(g_1,g_2,p)=(\ga_1,\gb(X), p_5(X))$.

Since $\sigma_{(d,0)}(\ga_1,\gb(X),p_5(X))=(\ga_1+\gb(X) d,\gb(X),p_5(X))$ for
$d\in A$, applying symplectic transvections of the type $\sigma_{(d,0)}$, say
$\Psi_2$, we may assume that if $\Psi_2(\ga_1,\gb(X),p_5(X))
=(\ga_2,\gb(X),p_5(X))$, then $\ga_2 \in R[X]$ and $\alpha_2= 1$ modulo
$s(X-1)R[X]$. Now, we write $\alpha_2$ as $\alpha_2(X)$.
Since $\beta(0)=1$, $(\alpha_2(X), \beta(X), p_5(X))\in \Um(R[X]^2 \perp
F_1,\langle ,\rangle)$.

Let $\beta(X) = 1 - X w$ and $\ga_2(X)=1+s(X-1)w'$ for some $w,w'\in
R[X]$. Then $s = sXw+ 
s\beta(X)$ and $\ga_2(X) = 1 + sXw' -(sXw+s\gb(X))w'$. Let $\ga_3(X)=
1 +sXw'(1-w)$. Then
$\sigma_{(sw',0)}(\ga_2(X),\gb(X),p_5(X))=(\ga_3(X),\gb(X),p_5(X))$ with 
$\alpha_3(X) = 1$ modulo
$sXR[X]$. 

Since $(\alpha_3(X),s)R[X] = R[X]$ and $ \beta(X)$ is monic, there
exists $c \in R$ such that $1- cs \in R \cap (\alpha_3(X), \beta(X)).$
Recall that $sM\subset F_1$. Therefore, writing $b = 1, b^\prime = 1-sc$ 
and applying (\cite{Bh-2}, Lemma 4.7),
there exists $\Psi_3 \in SL_2(R[X],(sX))\,E(R[X]^2 \perp M,
\langle,\rangle)$ such that 
$$\Psi_3 (\alpha_3(X), \beta(X), p_5(X)) =
(\alpha_3(b^\prime X), \beta(b^\prime X), p_5(b^\prime X)).$$ 
Since $
\alpha_3(X)= 1$ modulo $(sX)R[X]$, $\alpha_3(b^\prime X) =1$
modulo $(sb^\prime X)R[X]$. Moreover $b^\prime = 1- cs \in R \cap
(\alpha_3(b^\prime X), \beta(b^\prime X)).$ Therefore $[\alpha_3(b^\prime
X), \beta(b^\prime X)]$ is a unimodular row.

Let $\Psi_3 = {\Delta}^{-1} \Phi$, where $\Delta \in \SL_2(R[X], (sX))$
and $\Phi \in E(R[X]^2 \perp M, \langle,\rangle).$ Let 
$\Delta(\ga_3(b'X), \gb(b' X)) = (\ga_4(X),
\gb_1(X))$. Then $$\Phi(\ga_3(X), \gb(X), p_5(X)) =
(\ga_4(X), \gb_1(X), p_5(b' X)).$$ 
Since $\Delta \in
\SL_2(R[X], (sX))$, hence $ \ga_4(X) =1$ modulo $(sX)R[X]$ and
$[\alpha_4(X), \beta_1(X)]$ is a unimodular row. 

Write $\Gamma_2 = (\Phi\ot A)(\sigma_{(sw',0)}\ot A) \Psi_2 \Gamma_1$. Then
$\Gamma_2\in Esp(A^2\perp P,\langle,\rangle)$ and 
$\Gamma_2 (g_1,g_2,p)=(\ga_4(X),\gb_1(X),p_5(b'X))$ with
$[\alpha_4(X), \beta_1(X)]$ a unimodular row. 
Therefore, by
(\cite{Bh-2}, Lemma 4.1), there exists $\Phi_1 \in ESp(A^2 \perp P,
\langle,\rangle)$ such that
$$\Phi_1 (\alpha_4(X), \beta_1(X), p_5(b^\prime X))
= (\alpha_4(X), \beta_1(X), e_1).$$ 
Since $\langle e_1,f_1\rangle =
s$, $(\alpha_4(X), e_1)$ is an element of $\Um(A \oplus P).$ Therefore,
by (\cite{Bh-2}, Lemma 4.4), there exists $\Phi_2 \in ESp(A^2 \perp P,
\langle,\rangle)$ such that $\Phi_2 (\alpha_4(X), \beta_1(X), e_1) =
(1, 0, 0).$

Let $\Gamma = \Phi_2 \Phi_1 \Gamma_2$. Then $\Gamma (g_1,
g_2, p) = (1,0,0).$ Hence, the theorem is proved.  $\hfill \square$
\end{proof}

The proof of the following result follows from 
(\cite{Bh-2}, Lemma 5.2 and 5.4) and (\ref{88}).

\begin{theorem}\label{87}
Let $R$ be a ring of dimension $2$
and $A=R[X_1,\ldots,X_r,Y_1^{\pm 1},\ldots,Y_{r'}^{\pm 1}]$. Let
$P$ be a projective $A$-modules of rank $2$ with trivial determinant.
If $A^2$ is cancellative, then $P$ is cancellative.
\end{theorem}

\begin{proposition}
Let $R$ be a smooth affine domain of dimension $2$ over an
algebraically closed field $k$ of characteristic $0$. 
Let $A=R[X_1,\ldots,X_n,Y^{\pm
1}]$. Then $A^2$ is cancellative and hence every projective $A$-module
of rank $2$ with trivial determinant is cancellative (\ref{87}).
\end{proposition}

\begin{proof}
Let $P$ be a stably free $A$-module of rank $2$.  By (\ref{lindel}),
we may assume that $P\op A\iso A^3$. Since $A_{1+Yk[Y]} = \wt
R[X_1,\ldots,X_n]$, where $\wt R$ is a smooth affine domain over a
$C_1$ field $k(Y)$. Hence, by (\cite{Bh-2}, Theorem 5.5), $P\ot
A_{1+Yk[Y]}$ is free. There exists $h\in 1+Yk[Y]$ such that $P_h$ is
free. Patching $P$ and $A^2_h$, we get a projective
$R[X_1,\ldots,X_n,Y]=B$-module $Q$ of rank $2$ such that $Q_h
\iso P$. Since $(Q\op B)_Y$ is free,
$Q\op B$ is free. Applying (\cite{Bh-2}, Theorem 5.5), $Q$ is free and
hence $P$ is free. This proves that $A^2$ is cancellative.
  $\hfill \square$
\end{proof}
\vspace*{.1in}

\noindent{\bf Acknowledgments.}
I sincerely thank S.M. Bhatwadekar and
Raja Sridharan for suggesting the problem and
for useful discussion. In-fact it was Raja's idea that we should invert just 
special monic polynomials in (\ref{M2}) instead of all monic polynomial as 
was the case in \cite{Das}. 


{}

\end{document}